\documentclass[reqno,b5paper]{amsart}
\usepackage{amsfonts}
\usepackage{amssymb}
\usepackage{amsmath}
\usepackage{amsthm}
\usepackage{pgf,tikz}
\usepackage{capt-of}
\usepackage{syntonly}
\usepackage{pstricks,subfigure,epsfig}
\usepackage{amsmath,amsfonts,amssymb,amsthm,euscript,upgreek}
\usepackage{enumerate}
\usepackage{multirow}
\usepackage{appendix}

\theoremstyle{plain}
\newtheorem{thm}{Theorem}[section]
\newtheorem{lemma}{Lemma}[section]

\newtheorem{remark}{Remark}[section]
\newtheorem{cor}{Corollary}[section]
\newtheorem{defn}{Definition}[section]

\begin{document}

\title[The first law of cubology for the Rubik's Revenge] {The first law of cubology for the Rubik's Revenge}

\author[S.Bonzio]{Stefano Bonzio}
\address{Stefano Bonzio, University of Cagliari\\
Italy}
\email{stefano.bonzio@gmail.com}

\author[A.Loi]{Andrea Loi}
\address{Andrea Loi, University of Cagliari\\
Italy}
         \email{loi@unica.it}

 \author[L.Peruzzi]{Luisa Peruzzi}
\address{Luisa Peruzzi, University of Cagliari\\
Italy}
\email{luisa$\_$peruzzi@virgilio.it}

\thanks{
The second author was  supported by Prin 2010/11 -- Variet\`a reali e complesse: geometria, topologia e analisi armonica -- Italy and also by INdAM. GNSAGA - Gruppo Nazionale per le Strutture Algebriche, Geometriche e le loro Applicazioni.}

\subjclass[2000]{Primary: 05EXX; Secondary: 20BXX}
\keywords{Combinatorial puzzles; Rubik's Revenge}
\date{\today}

\begin{abstract}
We state and prove the first law of cubology of the Rubik's Revenge and provide necessary and sufficient conditions for a randomly assembled Rubik's Revenge to be solvable.
\end{abstract}

\maketitle




\section{Introduction} 

Erno Rubik, in 1974, invented the most famous and appreciated puzzle of all times that still goes under his name as Rubik's Cube. A few years later, in 1981, Peter Sebesteny, following Rubik's idea, invented his own cube, called the Rubik's Revenge, as it was meant to be a more difficult puzzle with respect to the original. 

The Rubik's Cube attracted the attention of many mathematicians (see, e.g. \cite{Bande82}, \cite{Joyner08}, \cite{Kosniowski81})  and in most cases they successfully gave a group theoretical analysis and solution to the puzzle. 

Any \emph{Cubemaster} knows that dismantling the cube and reassembling it randomly may cause in most of the cases the puzzle not to be solvable anymore. A question arises naturally: under which conditions a cube is solvable? The answer came a few years after the Cube was born. Indeed, Bandelow \cite{Bande82} has provided necessary and sufficient conditions for the solvability of the cube in a Theorem which he has christened  ``the first law of cubology'' (see Theorem \ref{3per3} below). As far as we know, the same question has not been answered for the Rubik's Revenge. Our aim is then providing an answer to this question. 

The paper is structured as follows: in the first section we present a review of the main interesting results, in our perspective, concerning Rubik's Cube. In the second section we go through the analysis of Rubik's Revenge and state Theorem \ref{risRevenge}, which  represents ``the first law of cubology'' for it, and prove some corollaries of our analysis. For example we provide necessary and sufficient conditions for a randomly assembled Rubik's Revenge to be solvable.
Our proof is based on the algebraic tools on the Rubik's Revenge developed in \cite{Larsen85}.  To  authors' best knowledge reference  \cite{Larsen85}  is the only place where a group theoretical approach to the Revenge is given (the reader is referred e.g.  to \cite{Adams82} and to several places in the web for the description of the instructions needed to  solve the  Rubik's Revenge).

\section{Review on the Rubik's Cube}

The Rubik's cube is composed by 26 small cubes, which we will refer to as ``cubies'' (as in \cite{Chen}). A quick look and one can notice that 8 are \textit{corner} cubies, i.e. cubies with 3 visible coloured faces, 12 are \textit{edge} cubies, with just 2 visible faces and the remaining 6 have one visible face: the \textit{center} cubies. 

The cube, obviously, has 6 faces, each of which can be moved clockwise or anticlockwise. Moving a face implies the movement of any of the cubies lying in the moved slice, with the exception of the center piece, occupying the same (spatial) position (\textit{centers} are fixed). 

Solving the cube means having every face of a unique colour: centers, being fixed, establish which colour the face shall have. For example, if one sees a face with the center coloured in white, when the cube is solved, the whole face will be of white colour. Of course, the same applies to all the other faces. 

Any \emph{Cubemaster} knows that if the cube is disassembled and then reassembled randomly, it may happen that it is not solvable anymore, as pieces shall be assembled following a precise pattern. On the other hand, mathematicians know that such problems can be studied using group theory \cite{Joyner08}, \cite{Larsen85}, \cite{Kosniowski81}, \cite{Signm77}, \cite{Signm82}. 

It is easily verified that the moves of the cube form a group, generated by the basic moves, generally referred to as $ R,L,F,B,U,D$ (as in \cite{Bande82} and \cite{Kosniowski81}). 

Corner and edge cubies can be moved as well as twisted, so they can change position (in space) as well as orientation. Hence, a natural way to describe a pattern is introducing permutations for position changes and orientation for twisting. In this way, a random pattern corresponds to a configuration, that can be described by a 4-tuple $ (\sigma , \tau , x , y ) $, as done in \cite{Bande82}. Permutations involving corner cubies are necessarily disjoint from the ones involving edges, as this is imposed by the construction of the cube itself. So $ \sigma $ refers to a permutation of corners, while $ \tau $ is a permutation on edges. Thus, in principle, $ \sigma\in S_{8} $, while $ \tau\in S_{12} $. When the cube is solved, clearly $ \sigma=id_{S_{8}} $ and $ \tau=id_{S_{12}} $.

Orientations can be characterized using vectors. As corners are eight and they have three visible faces, they may assume three possible different orientations, so the vector describing corners' orientation is $ x\in (\mathbb{Z}_{3})^{8} $; while edge cubies are twelve, but they have only two possible orientation, the vector is $ y\in (\mathbb{Z}_{2})^{12} $. 

Let us make clear how to calculate a random configuration of the Rubik's Cube. We assume the convention that we look at the Cube in order to have the white face on top and the red in front. Then we associate a number to the spatial position of each corner as well as of each edge. We assign a number from 1 to 8 to the position occupied by each corner \footnote{The idea is suggested by Bandelow \cite{Bande82} who uses the suggestive terminology of ``second skin'' for the spatial positions occupied by cubies.}. We number 1 the up-front-left corner and then associate numbers 2, 3, 4 just counting the others standing in the upper face clockwise. For corners standing in the down face, the down-front-left corner is assigned number 5 and the others take 6, 7 and 8 counting clockwise. \vspace{5pt}
\begin{center}
\includegraphics[scale=0.08]{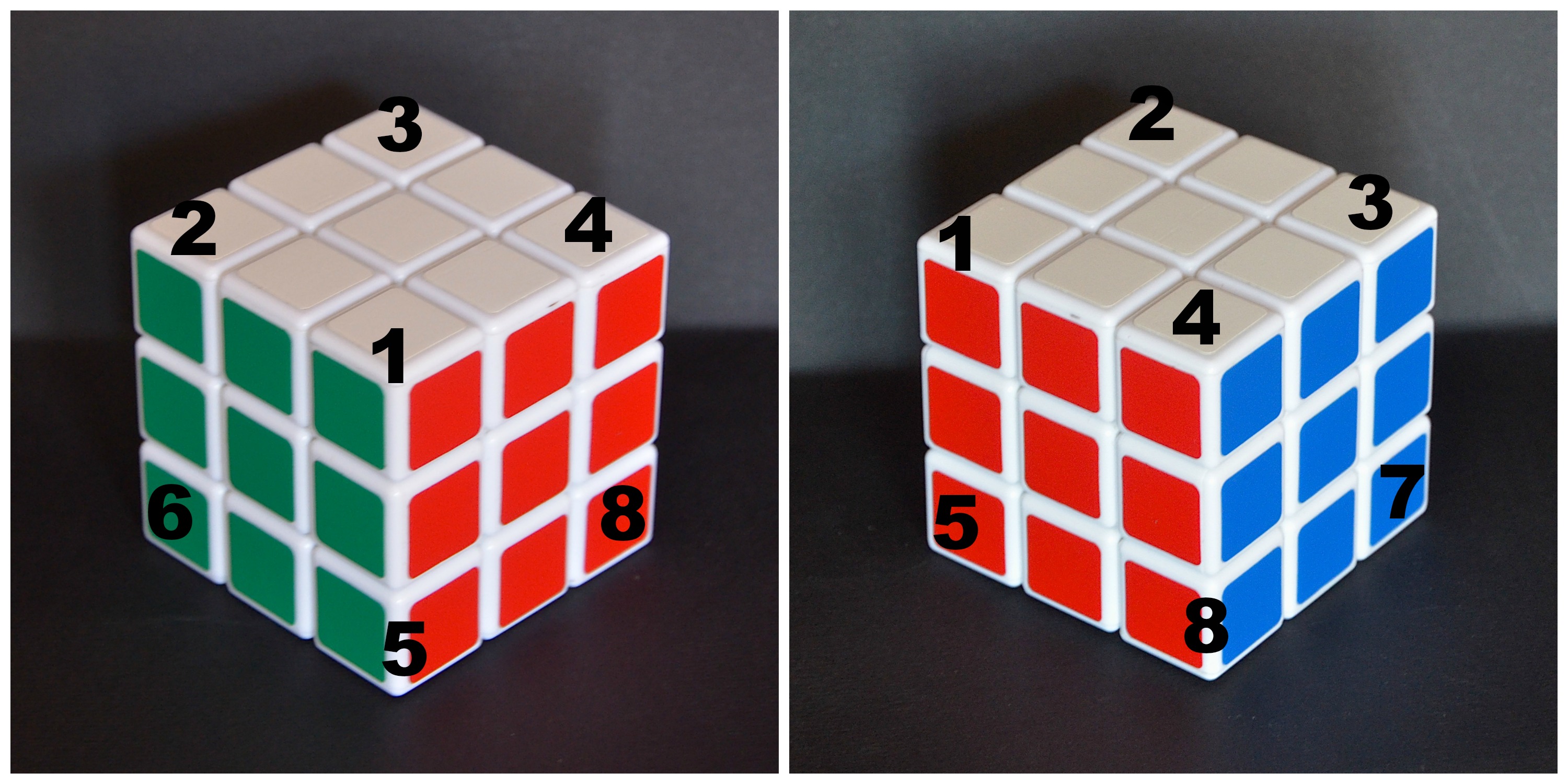}
\captionof{figure}{Enumeration of the spatial positions occupied by corners}
\end{center}

The same can be done with respect to edges: we assign numbers from 1 to 4 for the spatial positions of edges in the up face, starting from the front-up and counting then clockwise on the upper face. We number 5 the front-left position in the middle layer and then counting clockwise we give numbers 6, 7 and 8 to the others in the same layer. Finally we assign numbers from 9 to 12 to edge spatial positions in the down face, with 9 assigned to the front-down and the other counting clockwise. 
\begin{center}
\includegraphics[scale=0.10]{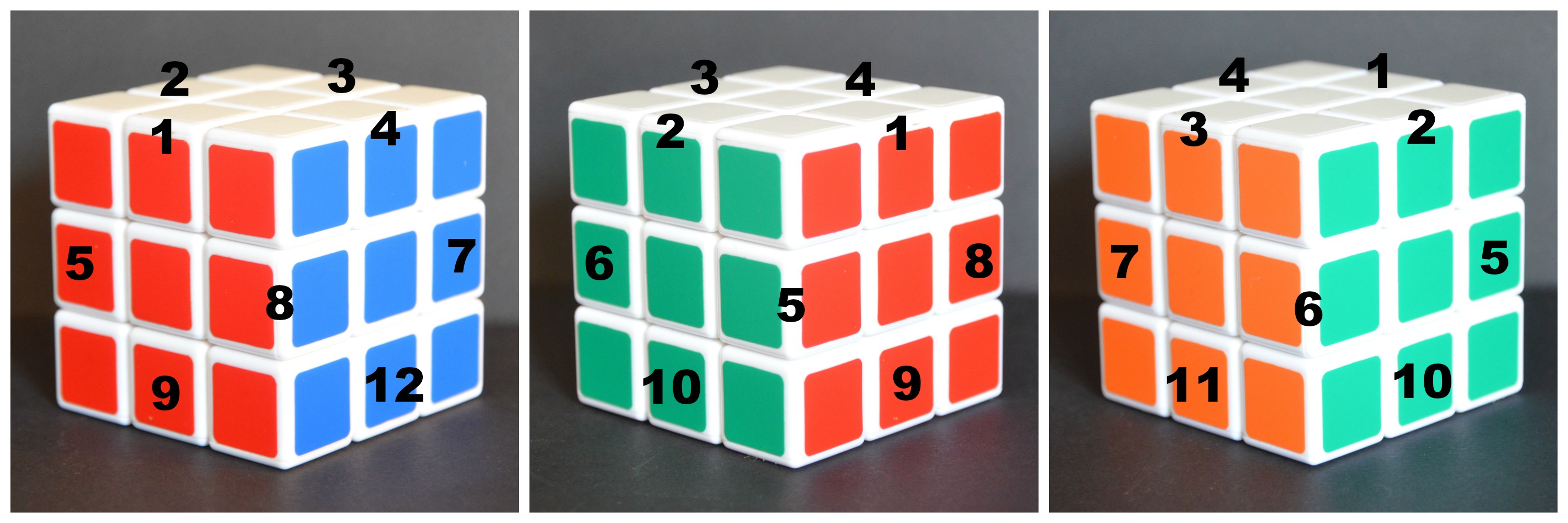}
\captionof{figure}{Enumeration of the spatial posistions for edges.}
\label{numeri edges}
\end{center}

Now let us see how to assign $ x_{i}\in\mathbb{Z}_{3} $ for each $ i\in\{1,...,8\} $. First of all, we decide that for corners having a white sticker, the latter is assigned with number 0 and the other stickers take number 1 and 2 moving clockwise on the cubie's faces, starting from the white one, as done in \cite{Chen}. Similarly, for corners having a yellow sticker, it takes number 0 and the others 1 and 2 counting clockwise. 
\vspace{5pt}
\begin{center}
\includegraphics[scale=0.08]{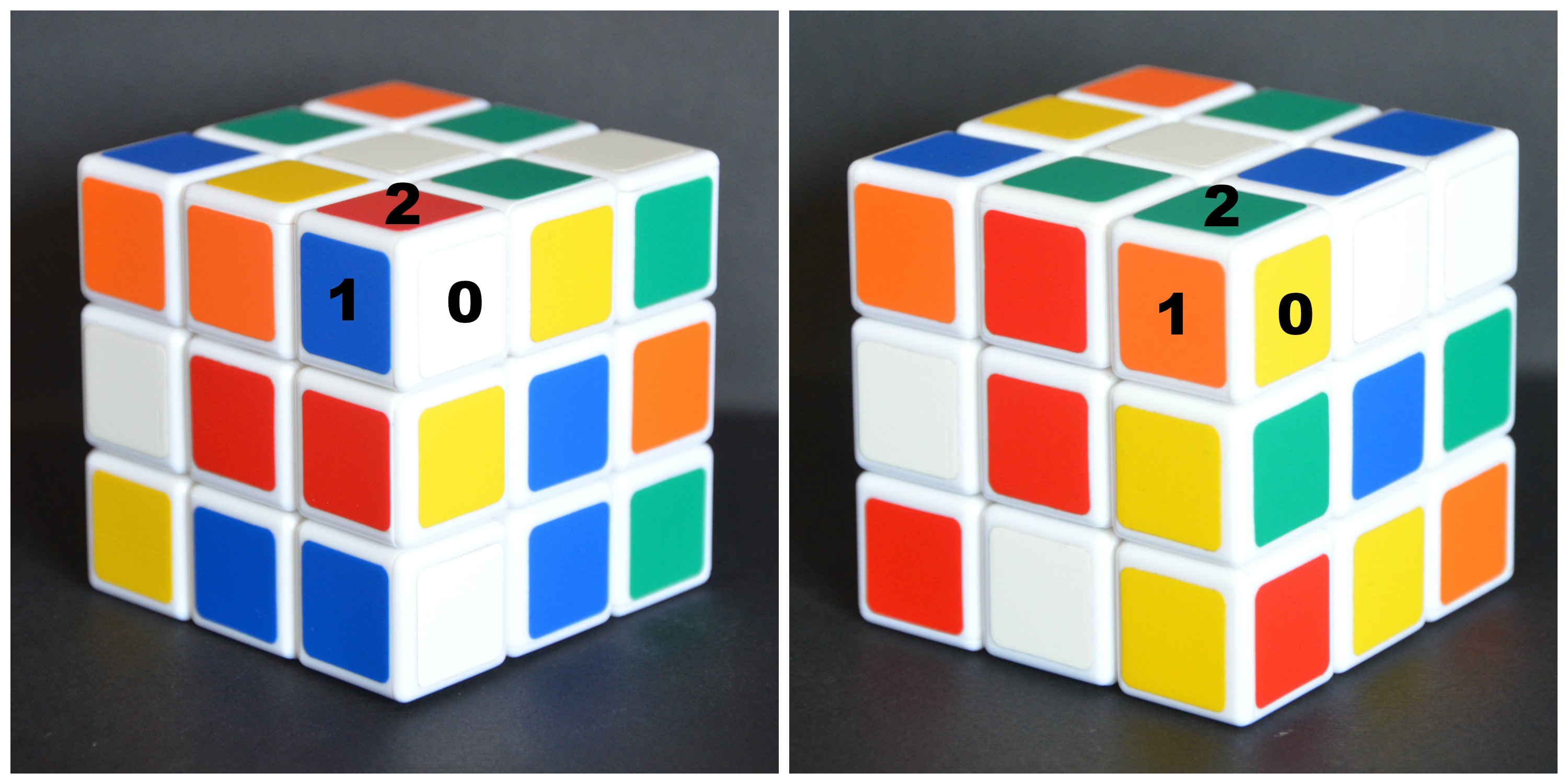}
\captionof{figure}{Example of assignation of numbers to stickers for corner cubies.}
\label{n°orientazione corner}
\end{center}

The convention for calculating components of vector $x$ for a random configuration is the following: we associate to the corner living in the i-th spatial position the orientation number $ x_{i}\in\mathbb{Z}_{3} $, defined as the number of the corner's sticker lying on the white or yellow face of the cube. Referring for example to the random configurations illustrated in both sides of Fig.\ref{n°orientazione corner}, we would get $ x_{4}=2 $, as the corner standing in position four has the stickers that takes 2 in the up face of the cube.

We proceed similarly for edges, i.e. we establish that for edges having a white or a yellow sticker, those ones take number 0 and the other stickers take 1 (examples in the left-hand side of Fig. \ref{n° orientazione edge}). For the remaining 4 edges, we decide that red and orange stickers take 0, while green and blue ones take 1 (examples in the right-hand side of Fig. \ref{n° orientazione edge}).  
\vspace{5pt}
\begin{center}
\includegraphics[scale=0.08]{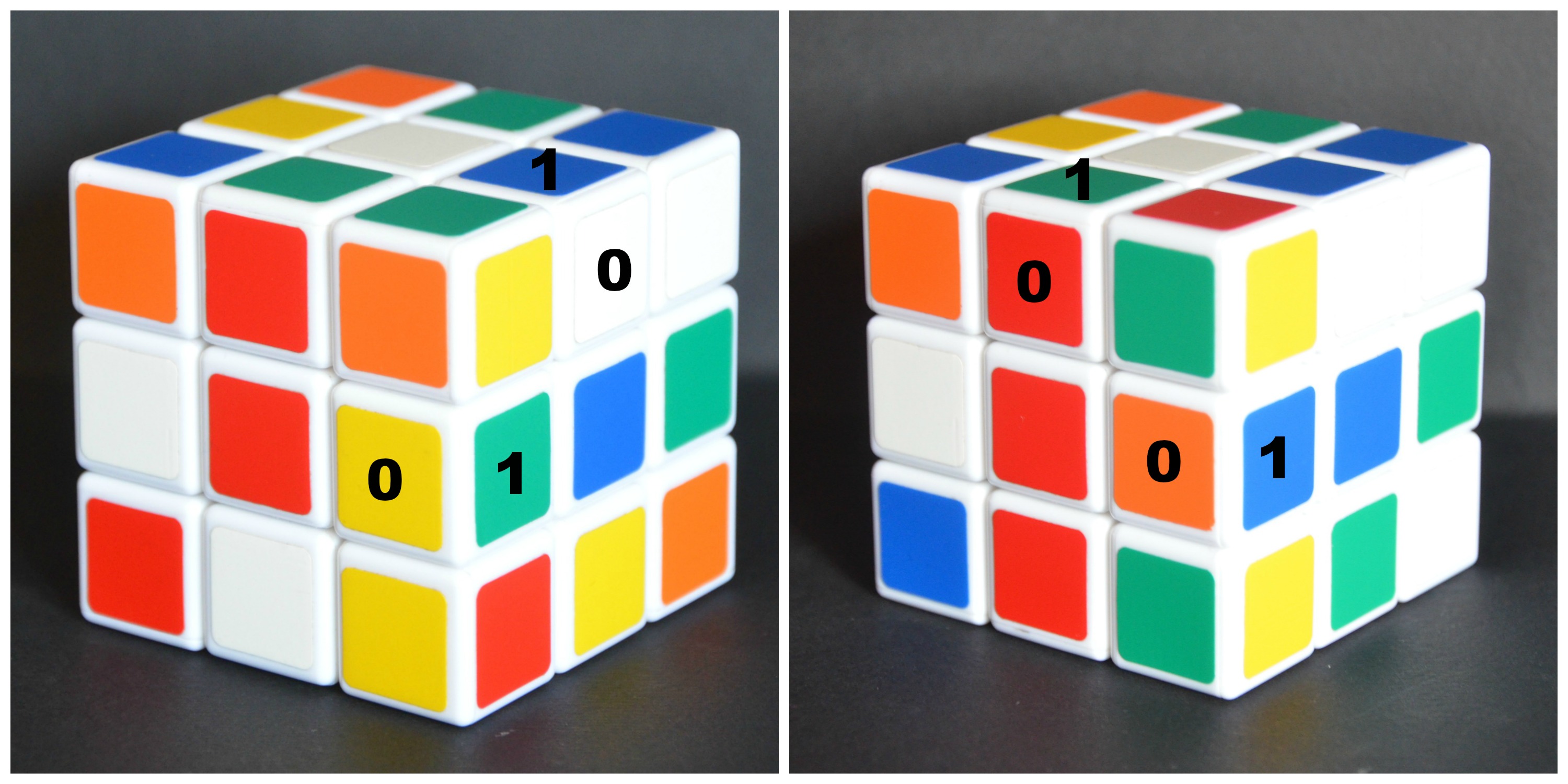}
\captionof{figure}{Assignation of numbers to stickers for edge cubies.}
\label{n° orientazione edge}
\end{center}
Determining $ y $ is done fixing four sides to look at the cube: up (white), down (yellow), front (red) and back (orange). We associate to the edge living in the i-th spatial position the orientation number $ y_i\in\mathbb{Z}_{2} $, defined as the number of the edge sticker lying respectively on the white, yellow, red or orange face of the cube.  For the random configuration illustrated in the left-hand side of Fig. \ref{n° orientazione edge}, we have $ y_4= 1 $, as the edge occupying position 4 has its blue sticker (taking number 1) lying in the upper face of the Cube; similarly $ y_8 =0 $ since the edge in position 8 has its yellow sticker (taking number 0) living in the front face of the Cube. Following the same principle for different stickers' colours and for the random configuration illustrated in the right-hand side of Fig. \ref{n° orientazione edge}, we have $ y_{1}=1 $ and $ y_{8}= 0 $. It is clear for the convention we have introduced that whenever the Cube is solved we have $ x_{i}= 0 $ for all $ i\in\{1,...,8\} $ and $ y_i=0 $ for all $ i\in\{1,...,12\} $.
   
We say that a configuration is \emph{valid} if one can reach the configuration  \\
$(id_{S_{8}}, id_{S_{12}}, 0 ,0)$, i.e. the configuration where the cube is solved, by a finite number of moves. \\
The ``first law of cubology'' \cite{Bande82} (Theorem 1, page 42) provides necessary and sufficient conditions for a configuration to be valid.  
\begin{thm}\label{3per3}
\textbf{\emph{(First law of cubology)}}
A configuration of the Rubik's Cube is valid if and only if
\begin{itemize}
\item[i)] $ sgn(\sigma) = sgn(\tau)$;
\item[ii)] $ \sum_{i}x_{i}\equiv 0 $(mod 3);
\item[iii)] $ \sum_{i}y_{i}\equiv 0 $(mod 2).
\end{itemize}
\end{thm}

From this theorem we get that the probability for  a  randomly assembled Rubik's cube to be solvable is \begin{large}$\frac{1}{12}$\end{large} and hence one gets:
\begin{cor}\emph{(\cite{Bande82} Theorem 2, page 44)}\label{configurationi 3per3}
The total number of possible patterns\footnote{ By possible patterns we mean the ones in the valid configuration.} is \begin{Large} $ \frac{8!\;\cdot\;3^{8}\;\cdot\;12!\;\cdot\;2^{12}}{12} $. \end{Large}
\end{cor}

Some relevant mathematical properties of the Rubik's Cube follow from the above theorem, see \cite{Bande82} and \cite{Joyner08} for details.

\section{Configurations of the Rubik's Revenge}
The Rubik's Revenge has been created a few years after the original Rubik's cube: every face is composed by four slices instead of three. 
The Rubik's Revenge is composed by 56 cubies: 8 are corner cubies exactly as in the original Cube, 24 are edge cubies and the remaining 24 are center cubies. At first glance, the big difference with the Rubik's Cube is that center pieces are not fixed anymore (clearly, also the number of edges is duplicated). As center cubies can be moved, they do not determine which colour  a face shall assume solving the Revenge. However, it is enough to choose a random corner to determine the colour that every face shall assume. Throughout the paper we mean the Revenge oriented so to have the white face on top and the red one in front. Hence the white-red-green corner, for example, shall occupy the up-front-left position in the solved Cube. In order to have that, after a quick look to the Revenge, we search the white-red-green corner and establish that position one is exactly located where such corner is living in, hence we rotate the whole cube so to have such a corner standing in the up-front-left position.  

As for the original Cube, the set of moves naturally inherits the structure of a group, denoted by $ \mathbf{M} $.
This group is generated by the twelve clockwise rotations of slices denoted by $ R,L,F,B,U,D, C_{R}, C_{F}, C_{U}, C_{L}, C_{B}, C_{D} $ where $ R,L,F,B,U,D $ are twists of the external slices, respectively, right, left, front, back, up and down face, while $ C_{F}, C_{R}, C_{U}, C_{L}, C_{B}, C_{D} $ are the twists of the central-front, central-right, central-up, central-left, central-back and central down slice respectively (some of which is illustrated in Fig. \ref{mosse di base}).
\vspace{5pt}
\begin{center}
\includegraphics[scale=0.12]{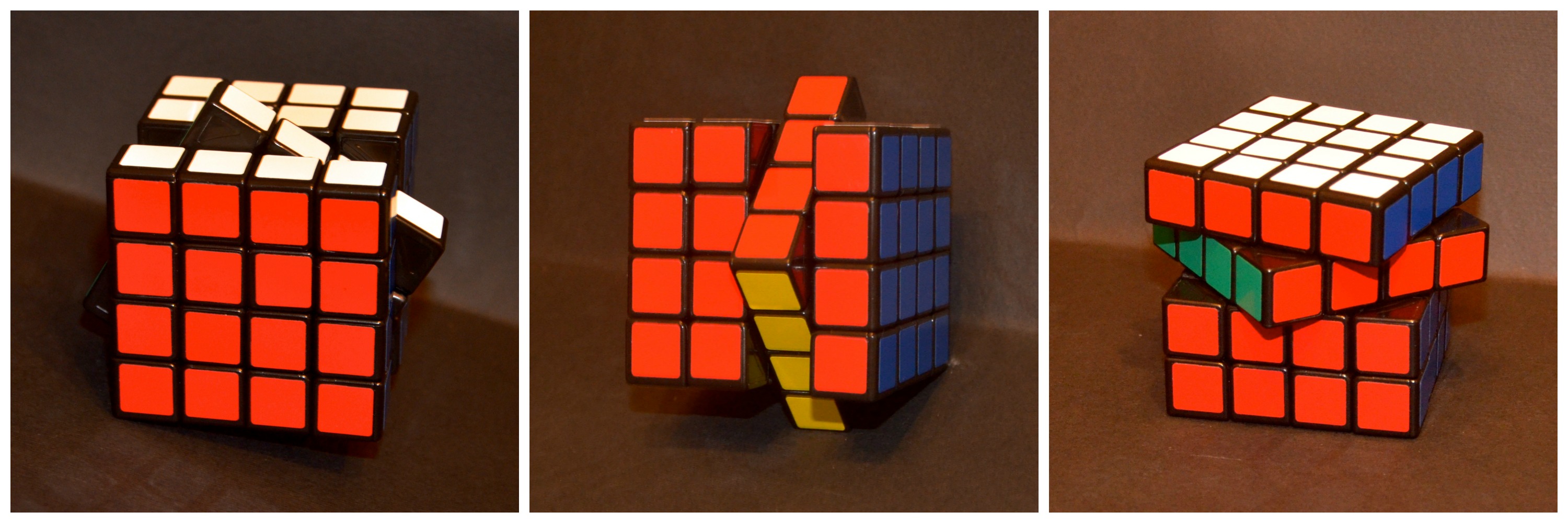}
\captionof{figure}{The moves $ C_{F}, C_{R} $ (left-hand side and central pictures) and the move $ C_{U}^{-1} $ (right-hand side) that is an anti-clockwise rotation of the central-up slice on the Revenge.} 
\label{mosse di base}
\end{center}

Any of those elements has order 4. 

The number of stickers one can find on the Revenge is equal to 96, so we may define an homomorphism
$$ \varphi :\mathbf{M}\longrightarrow \mathbf{S}_{96}, $$ 
\noindent
which sends a move $ m\in\mathbf{M} $ to a permutation $ \varphi(m)\in\mathbf{S}_{96} $ induced on the Revenge by the move $m$. The image $ \varphi (\mathbf{M}) \subset \mathbf{S}_{96} $ corresponds to those permutations in $ \mathbf{S}_{96} $ \textit{induced} by moves of the Revenge. 
\begin{remark}\label{inclusione stretta}\rm
\textrm{The inclusion $ \varphi (\mathbf{M}) \subset \mathbf{S}_{96}$ is (obviously) strict. For example, it may never happen that an element of 
$\mathbf{M}$ send a corner to the position occupied by a center or an edge.}    
\end{remark} 

We can then define the group of the Rubik's Revenge $ \mathbf{G}:=\varphi(\mathbf{M}) $. By the well known ``isomorphism theorem'' \cite{Hall76}, $ \mathbf{G}= \mathbf{M}/ker(\varphi) $, i.e. it is the group of the moves we get identifying all the combination of moves leading to the identical permutation.

Consider the subset of $\mathbf{S}_{96}$ corresponding to  permutations and/or orientation changes of corners, edges and center cubies. The set all of these permutations will be called
the  \textit{space of the configuration of the Revenge} and will be denoted by  $\mathcal{S}_{Conf}$. 

\begin{remark}\label{re:no flip}
\emph{It is known that a single edge cubie can not be flipped\footnote{ This is mathematically proved in \cite{Larsen85} (Theorem 2), and is also a physical constrain (see http://www.instructables.com/id/How-to-put-a-4x4-Rubiks-Cube-Together/, for a detailed description on the construction of the Rubik's Revenge).};
 however we may theoretically think to flip a single edge (and hence changing its orientation) by swapping its stickers.
 Therefore, unlike the Rubik's Cube, the cardinality of $\mathcal{S}_{Conf}$ is larger than the number of patterns we can get by dismantling and reassembling the cube.  
}
\end{remark}

Clearly $ \mathbf{G}\subset \mathcal{S}_{Conf}\subset \mathbf{S}_{96}$ and $ |\mathcal{S}_{Conf}| < 96! $. More precisely 
\begin{equation}\label{cardconf}
|\mathcal{S}_{Conf}|= (24!)^2\cdot 2^{24}\cdot 3^{8}\cdot 8!.
\end{equation}
The group $ \mathbf{G} $ acts on the left on $ \mathcal{S}_{Conf} $:
$$ \mathbf{G}\times\mathcal{S}_{Conf}\longrightarrow\mathcal{S}_{Conf} $$
$$ (g,s)\longmapsto g\cdot s $$
where $ \cdot $ represents the composition in $ \mathbf{S}_{96} $. This gives raise to a left action of $ \mathbf{M} $ on the space of configurations, by $ m\cdot s=g\cdot s $, where $ g=\varphi(m) $, and vice-versa. For this reason, from now on, we will not make any distinction between the two actions on $ \mathcal{S}_{Conf} $. Notice that  the action of $\mathbf{G}$ on $ \mathcal{S}_{Conf} $ is \textit{free} (in contrast with that of $\mathbf{M}$), i.e.  if $g\cdot s = s$ then $g=id $. Hence this action yields a bijection between the group $ \mathbf{G} $ and the orbit $ \mathbf{G}\cdot s=\{g\cdot s\ | \ g\in\mathbf{G}\} $ of an arbitrary $s\in\mathcal{S}_{Conf} $, obtained by sending $g\in\mathbf{G}$ into $g\cdot s\in\mathcal{S}_{Conf} $.

\begin{remark}
\emph{It is easily seen that  the space of configuration $\mathcal{S}_{Conf} $ is a subgroup of $\mathbf{S}_{96}$ containing $\mathbf{G}$ as a subgroup. Then the left action $g\cdot s$, $g\in \mathbf{G}$ and $s\in \mathcal{S}_{Conf}$, can be also seen as the multiplication in $\mathcal{S}_{Conf}$ and 
the orbit  $\mathbf{G}\cdot s$ of $s\in \mathcal{S}_{Conf}$ is nothing but the right coset of $\mathbf{G}$ in  $\mathcal{S}_{Conf}$ 
with respect to $s$.}
\end{remark}

In order to characterize mathematically the notion of configuration we have to label center and edge cubies. Indeed, corners are univocally identified by the colours of their faces, but ambiguity may arise concerning edges and centers. 

In order to describe positions of each edge (or center) by permutations, we have to label all of them. A number between 1 and 24 is the label for center cubies.
Once all centers have been marked, the position of each of them in a random pattern can be described by a permutation $ \rho\in S_{24}$.

Regarding corners, things work like in the Rubik's cube, so a permutation $ \sigma\in S_{8} $ describes their positions. 

The twenty-four edge cubies can be devided in twelve pairs, namely those ones with the same colour. The two members of a pair are labelled with different letters: \textit{a} and \textit{b}, respectively. This is enough to provide a description of edges' positions by using a permutation $ \tau\in S_{24} $. We refer to an edge labelled with \textit{a} (respectively \textit{b}) as an edge of \textit{type a} (respectively \textit{type b}). Obviously the type is not dependent on the position the edge is lying in. 

Describing orientations of corners can be achieved by a vector $ x\in(\mathbb{Z}_{3})^{8}$ in the same way described in the previous section.
\begin{remark}
\emph{Due to the convention introduced above that the white-red-green corner is always set in the up-front-left position, i.e. in position 1, it will always happen that $ x_1= 0 $.}
\end{remark}
 In order to introduce such a vector for edges' orientation, we have to describe the spatial positions for edges. We proceed as done for the Rubik's cube (see Fig. \ref{numeri edges}), by using only twelve numbers (instead of 24) and the label \textit{a} and \textit{b}. 
\vspace{5pt}
\begin{center}
\includegraphics[scale=0.115]{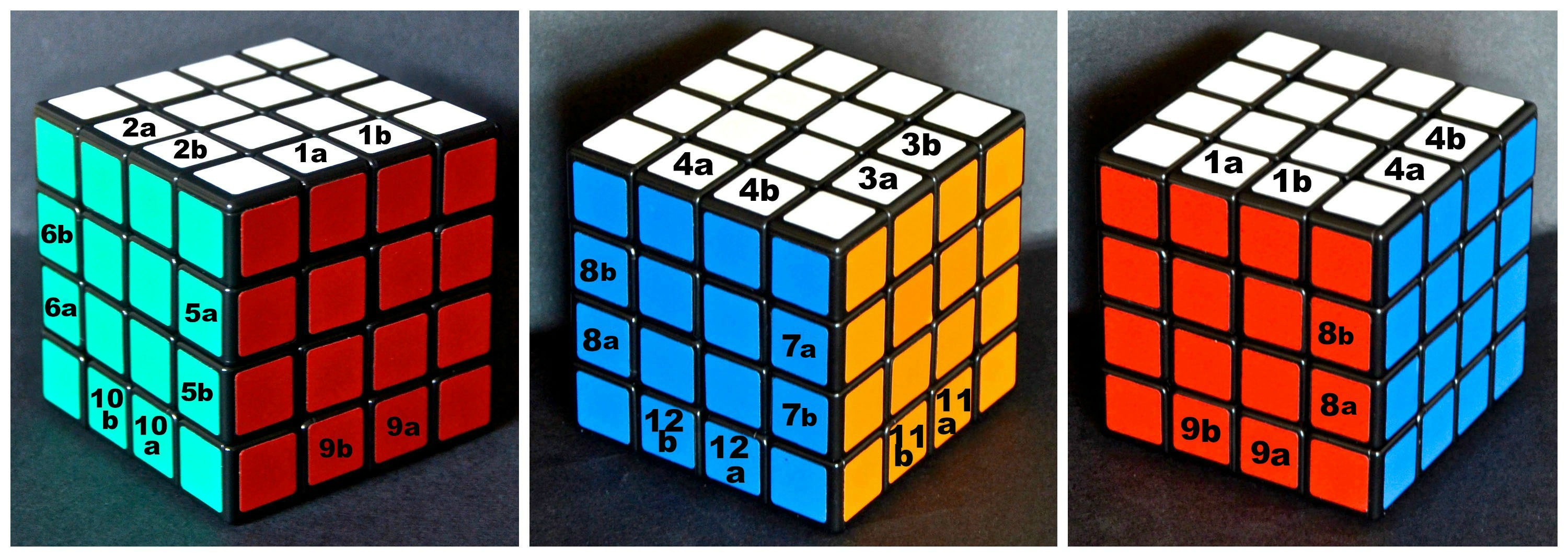}
\captionof{figure}{Schema of the assignation of numbers to the spatial position of edge cubies, by using labels \textit{a,b} (and twelve numbers).}
\label{assegnazioni a,b}
\end{center}

Concerning edges orientation, we can do the same as for the Rubik's Cube, shifting to two 12-tuple $ y_a=(y_{1_a},y_{2_a},...,y_{{12}_a}) $, with $ y_{i_a}\in\mathbb{Z}_{2} $ and $ y_b=(y_{1_b},y_{2_b}, $ \\
$...,y_{{12}_b}) $, with $ y_{i_b}\in\mathbb{Z}_{2} $ (edges are twenty-four, divided in pairs a and b).\\

It follows that the space of configurations $\mathcal{S}_{Conf}$  is in bijection with  the set of $5$-tuples $ (\sigma , \tau , \rho , x , y)$, where $ \sigma \in S_{8} $, $ \tau\in S_{24} $ , $ \rho\in S_{24} $, while $ x\in (\mathbb{Z}_{3})^{8} $ and $ y\in(\mathbb{Z}_{2})^{24} $.  
From now on we identify $\mathcal{S}_{Conf}$ with such $5$-tuples. The $5$-tuple $(id_{S_{8}}, id_{S_{24}}, id_{S_{24}}, 0 ,0)$
will be called the \textit{initial configuration}.

At the beginning of the section we claimed that the group of the moves \textbf{M} is generated by twelve elements. In fact, having introduced the formal notion of configuration, we may notice that nine generators are enough, as three moves involving central slices can be constructed as compositions of other basic moves. Notice, for example that $ C_{L}:=L^{-1}C_{R}R $ (center-left) and the same applies to other central moves that we will refer to as $ C_{B} $ (center-back) and $ C_{D} $ (center-down).

\begin{defn}\label{conf Revenge valida}
A configuration of the  Rubik's Revenge is \emph{valid} when it is in the orbit of the initial configuration under the action of $ \mathbf{G}$.
\end{defn}

We now present some basic facts concerning orbits, needed in the proof of Theorem \ref{risRevenge}.
\begin{lemma}\label{sg(sigma)=sg(tau)}
If two configurations $ (\sigma , \tau , \rho , x , y) $ and $ (\sigma' , \tau' , \rho' , x' , y') $ are in the same orbit then $ sgn(\sigma)sgn(\rho) = sgn(\sigma ')sgn(\rho') $. 
\proof
If $ (\sigma, \tau, \rho, x,y) $ and $ (\sigma ', \tau ', \rho', x',y') $ are in the same orbit, then $ (\sigma', \tau', \rho', x', $ \\
$y')= g\cdot(\sigma, \tau, \rho, x,y) $, for some $ g\in\mathbf{G}$. Hence it is enough to show that basic moves $ R,L,F,B,U,D $,\\ $ C_{R}, C_{F}, C_{U}$ preserve condition $ sgn(\sigma)sgn(\rho)=sgn(\sigma ')sgn(\rho ') $.

The action of $ g $ on corners is disjoint from the action of centers. It is easy to notice that any move among $ \{ R,L,F,B,U,D \} $ consists of a 4-cycle on both corners and centers, hence $ sgn(\sigma')=-sgn(\sigma) $ and $ sgn(\rho')=-sgn(\rho ) $, hence $ sgn(\sigma)sgn(\rho)=sgn(\sigma ')sgn(\rho ') $. On the other hand, moves $ C_{R}, C_{F}, C_{U} $ are identities on corners and consist of two 4-cycles on centers, implying that $ sgn(\sigma)=sgn(\sigma') $ and $ sgn(\rho)=sgn(\rho') $, hence $ sgn(\sigma)sgn(\rho)=sgn(\sigma ')sgn(\rho ')$ also in this case.    
\endproof
\end{lemma} 

\begin{lemma}\label{orientaz. x Revenge}
If $ (\sigma, \tau, \rho, x,y) $ and $ (\sigma ', \tau ', \rho', x',y') $ are configurations in the same orbit then $ \sum x'_{i}\equiv\sum x_{i} $ (mod 3).
\proof
The proof works exactly as for the anologuos property valid for the Rubik's Cube \cite{Bande82}.
\endproof
\end{lemma}

Before stating the main result of this paper, we make some considerations concerning edge cubies. We are aware of the fact that in a random configuration, an edge of type \textit{a} (resp. \textit{b}) can occupy either an  a-position or a b-position, as sketched for example in Fig. \ref{Revenge casuale}. Hence, using the information encoded in $ \tau\in\mathbf{S}_{24} $ we may associate to any edge a number $ i_{t,s} $, with $ t,s\in\{a,b\} $, where $ i_t $ indicates the spatial position, while $ s $ refers to the type of the edge. There are always orientation numbers associated to any edge $ i_{t,s} $ which will be $ y_{i_{t,s}}:= y_{i_t} $. 
  
\begin{center}
\includegraphics[scale=0.08]{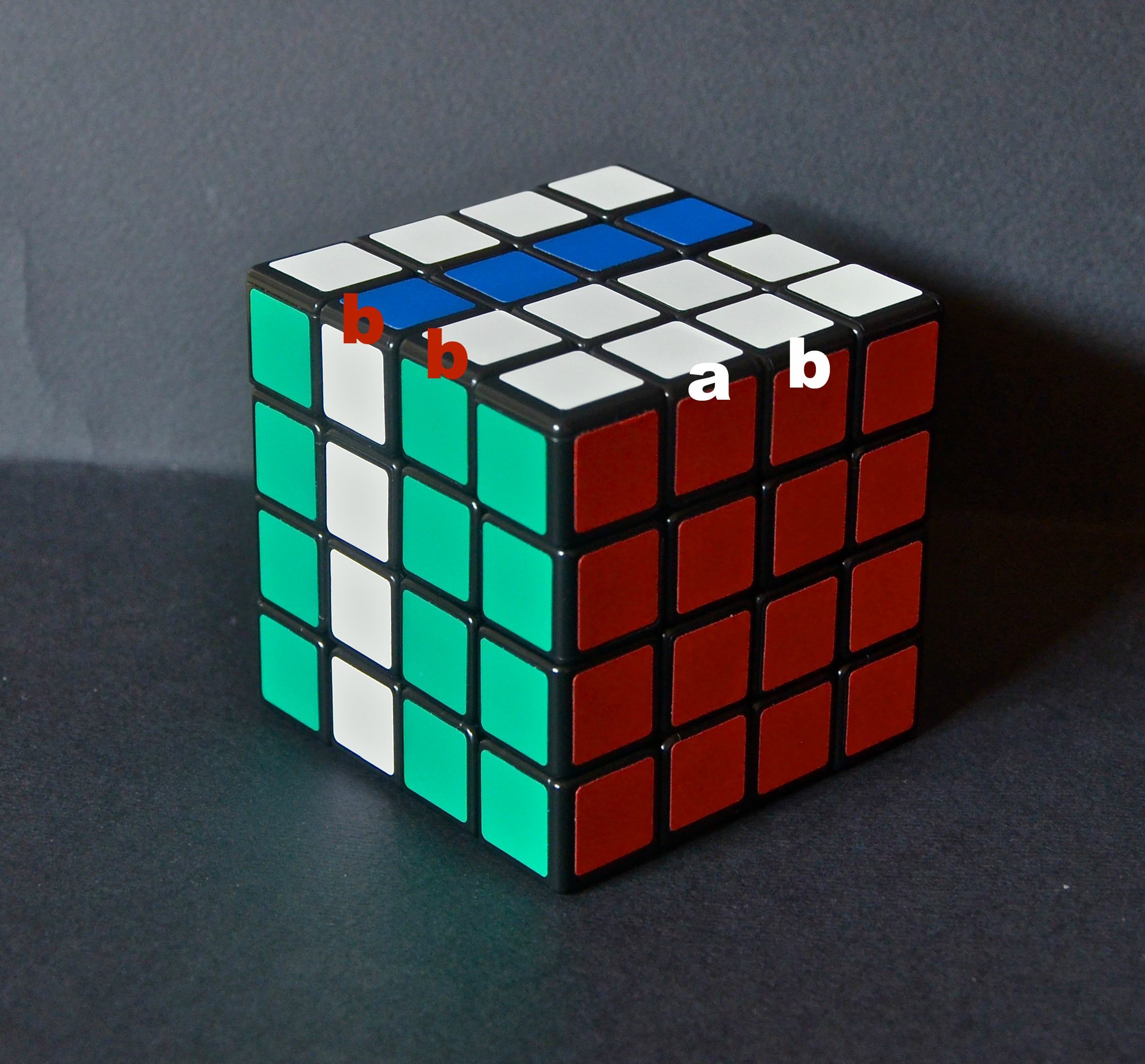}

\captionof{figure}{Letters a,b represents label associated to edges. In the above configuration we notice an edge of type \textit{a} in position $ 1_a $ and an edge of type \textit{b} in position $ 1_b $, hence the associated numbers will read $ 1_{a,a} $, $ 1_{b,b} $. In position 2 we find two edges of type \textit{b}, hence $2_{b,b}$ and $2_{a,b} $. Notice that the two edges in position 2 have two different orientations, namely $y_{2_{a,b}}=1$, $y_{2_{b,b}}=0$.}
\label{Revenge casuale}
\end{center}

We can now give the conditions for a configuration to be valid: this is actually the ``first law of cubology'' for the Rubik's Revenge:
\begin{thm}\label{risRevenge}
A configuration $ (\sigma , \tau , \rho , x , y) $ is valid if and only if 
\begin{enumerate}
\item[\emph{1.}] $ sgn(\sigma ) = sgn(\rho ) $ 
\item[\emph{2.}] $ \sum_{i}x_{i}\equiv 0 $(mod 3)
\item[\emph{3.}] $ y_{i_{t,s}}=1-\delta_{t,s} $, $ \forall i=1,...,12,$
\end{enumerate}
where  $\delta_{a,a}=\delta_{b,b}=1$ and  $\delta_{a,b}=\delta_{b,a}=0$.
\end{thm}   

Next section is dedicated to the proof of this theorem. Our proof will not be constructive, i.e. we do not show the moves actually needed to solve the cube, as we will use some group-theoretical results. Here we present some corollaries.  
\begin{cor}\label{ordine di G}
The order of $ \mathbf{G} $ is \begin{Large} $ \frac{(24!)^{2}\cdot\;8!\;\cdot\;3^{7}}{2} $ \end{Large}.
\proof
Since the action of $ \mathbf{G} $ is free, then $ |\mathbf{G}|=|\mathbf{G}\cdot s| $ for all $ s\in\mathcal{S}_{Conf.} $. It follows that $ |\mathbf{G}|=\frac{|\mathcal{S}_{Conf.}|}{N} $, where $ N$ is the number of orbits. Theorem \ref{risRevenge} yields $N=2\cdot 3\cdot 2^{24}$  and the results follows by (\ref{cardconf}).
\endproof
\end{cor}

The corollary above agrees with result presented in the last section of \cite{Larsen85}.\\

In order to study the solvability of the Revenge we give the following:

\begin{defn}\label{Revenge solvable}
A \emph {randomly assembled} Rubik's Revenge is a pattern of the Revenge obtained by any permutation and/or orientation change of corners, edges\footnote{ Here we allow also edge flips (see Remark \ref{re:no flip}).} and center cubies of the solved Revenge. 
\end{defn}

\begin{cor}\label{Sticker}
The probability that a randomly assembled Rubik's Revenge is solvable is \begin{large} $\frac{1}{2^{12}\cdot 3}$ \end{large}.
\proof
In a randomly assembled Rubik's Revenge center and edge cubies are not labelled, thus centers can always be moved so to have condition 1 in Theorem \ref{risRevenge} satisfied. The 24 equations in condition 3 are  reduced to 12: this can be obtained by assigning a label $a$ or $b$ to each edge in a pair, depending on its orientation, in such a way that $y_{i_{t,s}}=1-\delta_{t,s}$.
\endproof
\end{cor}

We have already mentioned (see footnote $3$)  the fact that the Revenge sold on  the market  is  different from 
the Revenge we consider in this paper.
In fact, the most relevant feature of the Revenge sold on the market is that any member of a pair of edges  of the same colours is different from its companion. This has the physical effect that it is impossible to assemble the cube putting an edge of type $a$ (respectively type $b$) in a '$b$-position' (resplectively $a$-position), without changing the orientation of both edges in a pair.  This yields that condition 3 in Theorem \ref{risRevenge}   can be always achieved, due to the internal mechanism of the Revenge.  Thus (surprisingly enough)  we get:
\begin{cor}\label{Monkey}
The probability that a randomly assembled Rubik's Revenge  sold on the market  is solvable is \begin{large} $\frac{1}{3}$ \end{large}.
\end{cor}

\section{Proof of Theorem \ref{risRevenge}}

We first introduce some already known results \cite{Larsen85} on the structure of $ \mathbf{G} $ which are of fundamental importance for proving Theorem \ref{risRevenge}.  

In order to prove the theorem we describe some significant subgroups of $ \mathbf{G}$, denoted by
$\mathbf{C}$, $\mathbf{Z}$ and $\mathbf{E}$ respectively. $\mathbf{C}$ is  the subgroup of $ \mathbf{G} $ which permutes corner cubies (no matter the action on orientation), and act as the identity on other pieces. $ \mathbf{Z} $ is the subgroup permuting centers only, leaving corners and edges fixed; $ \mathbf{E}$ permutes edges only.\footnote{Since each may assume three different orientations, it is known \cite{Signm82}, \cite{Larsen85} that the subgroup of corners corresponds to the wreath product $ \mathbf{H} = \mathbf{S}_{8} \bigotimes_{Wr} \mathbb{Z}_{3} $. However, we aim at describing the quotient subgroup $ \mathbf{C}= \mathbf{H}/\mathbf{T} $, where $ \mathbf{T} $ is the (normal) subgroup consisting of all al possible twists.}  
The proofs of these results can be found in \cite{Larsen85}. However they are included in the present paper for sake of completeness. Notice that our proofs are essentially the same, but in terms of the action of $ \mathbf{G} $ on $ \mathcal{S}_{Conf.} $.
\begin{thm}\label{C=A8}
$ \mathbf{C}\cong \mathcal{A}_{8} $, the alternating group of even permutation.  
\end{thm}
\noindent
We report to \cite{Signm82} for the proof of the above theorem. \\
In the next theorem we make use of the commutator, formally for $ m,n\in\mathbf{G} $, $ [m,n]=m\cdot n\cdot m^{-1}\cdot n^{-1} $. \\
Center cubies are 24, hence necessarily $ \mathbf{Z}\leqslant S_{24} $. 
\begin{thm}\label{Z=A24}
$ \mathbf{Z}\cong\mathcal{A}_{24} $, the alternating group of even permutation. 
\proof
We first show that $ \mathcal{A}_{24} \leqslant \mathbf{Z} $. The move 
\begin{equation}\label{mossa z}
z= [[C_{F},C_{D}],U^{-1}]
\end{equation}
is a 3-cycle on center and an identity on edges and corners. In fact, it is easy to check that the action of $ z $ on the initial configuration gives: $ z\cdot (id_{S_{8}}, id_{S_{24}}, id_{S_{24}}, 0, $ \\
$ 0)= (id_{S_{8}}, id_{S_{24}} , \rho_{1}, 0, 0) $, where $ \rho_{1} $ is 3-cycle. 

Observe that any three target centers can be moved to the positions permuted by $ z $ by a certain $ g\in\mathbf{G} $. Such a $ g $ admits an inverse $ g^{-1}\in  \mathbf{G} $, hence by $ g\cdot z\cdot g^{-1} $ we may cycle any center cubies. As $ \mathcal{A}_{24} $ is generated by any 3-cycle on a set of twenty-four elements, we have the desired inclusion. 

For $ \mathbf{Z}\leqslant\mathcal{A}_{24} $, we show that any odd permutation involving centers permutes necessarily also corners or edges, hence it cannot be in $ \mathbf{Z} $. Indeed, suppose that there exist $ \alpha\in\mathbf{Z} $ s.t. $ sgn(\alpha )=-1 $. Such an $ \alpha $ shall be obtained as a sequence of basic moves. Without loss of generality we can assume that $ \alpha $ is a sequence of $ L,R,U,D,F,B $, since the moves $ C_{R}, C_{F}, C_{U} $ consist of an even permutation on centers. On the other hand, $ L,R,U,D,F,B $ induces a 4-cycle on centers. However, all those moves must have permuted corners too, thus there exist a $ \beta=(\beta_{1},\beta_{2})\in \mathbf{S}_{8}\times \mathbf{E} $, s.t. $ sgn(\beta)= -1 $. Hence $ \beta_{1}\neq id_{S_8} $, implying that $ \alpha\not\in\mathbf{Z} $, the desired contradiction. 
\endproof 
\end{thm}
\noindent
Now we consider the subgroup $ \mathbf{E} $ of moves involving edges only. Edges are 24 each of which can assume two different orientation, however \textit{no single edge can be flipped} \cite{Larsen85}. The direct consequence of this fact is that $ \mathbf{E}\leqslant S_{24} $. Actually we can prove 
\begin{thm}\label{E=S24}
$ \mathbf{E}\cong \mathbf{S}_{24} $
\proof
We first show that $ \mathcal{A}_{24} \leqslant \mathbf{E} $. Indeed the move 
\begin{equation}\label{mossa e}
 e=[C^{-1}_{L},[L,U^{-1}]] 
 \end{equation}
is of a 3-cycle on edges. As done for centers, one can bring any target edge in the positions switched by $ e $ using an element of $ g\in\mathbf{G} $ and then solving the mess created by $ g^{-1} $. In this way, one obtains any 3-cycles in $ \mathbf{E} $, proving the inclusion. \\
We are left with proving that there is at least an odd permutation in $ \mathbf{E} $, which implies necessarily that $ \mathbf{E}\cong S_{24} $. \\
Consider the move $ C_{R} $: it gives raise to an even permutation on centers (two 4-cycles) and an odd one on edges (one 4-cycle). As by Theorem \ref{Z=A24} $ \mathbf{Z}\cong\mathcal{A}_{24} $, one can find $ z_{0}\in\mathbf{Z} $ such that $ z_{0}\cdot\varphi(C_{R}) $ acts as a 4-cycles of edges only. Hence $ z_{0}\cdot\varphi(C_{R})\in\mathbf{E} $ and $ sgn(z_{0}\cdot\varphi(C_{R}))=-1 $.  
\endproof   
\end{thm}

\noindent
\textbf{Proof of Theorem \ref{conf Revenge valida}} \\
\\
$ (\Rightarrow) $ Assuming $ (\sigma , \tau , \rho , x , y) $ is valid means that it is in the orbit of the initial configuration $ (id_{_{S_8}}, id_{_{S_{24}}}, id_{_{S_{24}}},  0, 0) $. \\
\\
\noindent
1. By Lemma \ref{sg(sigma)=sg(tau)} we have that $ sgn(\sigma)sgn(\rho)=sgn(id_{_{S_8}})sgn(id_{_{S_{24}}})=1 $, hence $ sgn(\sigma)= sgn(\rho) $, for their product must be equal to 1. \\
\\
2. $ \sum_{i}x_{i}\equiv 0 $(mod 3) follows trivially from Lemma \ref{orientaz. x Revenge} and the fact that $ (\sigma , \tau , \rho , x , y) $ is valid.  \\
\\
3. In the initial configuration $ (id_{S_{8}}, id_{S_{24}}, id_{S_{24}}, 0 ,0) $ it holds $ y_{i_t}= 0 $ for all $ i\in\{1,...,12\}$ and  $\delta_{a,a} =\delta_{b,b}= 1 $, hence $ y_{i_{t,s}}=1-\delta_{t,s}= 0 $. 

As $ (\sigma, \tau, \rho, x,y) $ is in the orbit of the initial configuration, it is obtained by a sequence of basic moves, thus we need to check that those moves preserve condition $ y_{i_{t,s}}= 1 - \delta_{t,s} $. 

We consider moves splitted in two sets: $ M_{1}=\{R, U, D, L\} $ and $ M_{2}=\{F, B, C_{R}, C_{F}, C_{U} \} $; hence we have two possibilities: we may assume a move, say $ m $, either $m\in M_{1}$ or $m\in M_{2}$. \\
Assume $m\in M_{1}$. Recall that for the convention we have introduced about the assignation of orientation numbers to edges, $ m $ does not change edge cubies' orientation, so we get $ y_{i_{t,s}}= 0 $ for all $ i\in\{1,...,12\} $. Furthermore $ m $ acts on a configuration moving edges occupying an a-position in edges in a-position and the same holds for b-positions and hence $\delta_{t,s} = 1$. 

Let now $m\in M_2 $. $ m $ changes orientations of some edges (the ones that it is actually permuting): more precisely it gives raise to a cycle of four edges or to two cycles of four edges each. Let $ i_{t,s} $ be one of those edges, then $ y_{i_{t,s}}= 1 $ and $ \delta_{t,s}=0 $ since a-positions and b-positions are swapped by $ m $.

$ (\Leftarrow ) $ We have to show that once conditions 1, 2 and 3 are satisfied we are always able to solve the cube. 
In the random configuration $ (\sigma, \tau, \rho , x,y) $ we can check (just by watching the Revenge) whether $ \rho $ is even or odd. If $ sgn(\rho)= -1 $, it is enough to apply one among $ \{ R,L,U,D,F,B\} $ to get $ sgn(\rho) = +1 $. 
If $sgn(\rho)= +1 $, then $ \rho\in \mathcal{A}_{24} $, hence, by Theorem \ref{Z=A24},  $\rho\in\mathbf{Z}\cong\mathcal{A}_{24} $, there exists $ z_{1}\in\mathbf{Z} $ s.t. $ z_{1}\cdot (\sigma, \tau, \rho , x,y)= (\sigma, \tau, id_{S_{24}}, x,y ) $. \\
By condition 1. $ sgn(\sigma) = sgn (\rho)=sgn(id) $, hence in $ (\sigma, \tau, id_{S_{24}} , x,y) $, $ sgn(\sigma) = +1 $. Then $ \sigma\in\mathcal{A}_{8} $. By Theorem \ref{C=A8}, $ \mathbf{C}\cong\mathcal{A}_{8} $, thus there exists $ c\in\mathbf{C} $ such that  $ c\cdot (\sigma, \tau, id_{S_{24}} , x,y) = (id_{S_{8}}, \tau, id_{S_{24}} , x,y) $. 

Now, we proceed setting edges in their correct positions; $ \tau $ is a permutation of 24 elements, however Theorem \ref{E=S24} $ \mathbf{E}\cong\mathbf{S}_{24} $, hence there is an $ e_{1}\in\mathbf{E} $ such that \\$ e_{1}\cdot (id_{S_{8}}, \tau, id_{S_{24}} , x,y)= (id_{S_{8}}, id_{S_{24}}, id_{S_{24}} , x,y) $. Condition (3)  implies that, as all edge cubies are correctly positioned, then they are also correctly oriented, thus $ y=0 $, so $ (id_{S_{8}}, id_{S_{24}} , id_{S_{24}} , x,y)= (id_{S_{8}}, id_{S_{24}}, id_{S_{24}} , x,0) $. 

Now the labelled Revenge has been reduced to Rubik's Cube, as any pair of edge can be seen as an unique big edge. So we have actually reduced the Rubik's Revenge to a Rubik's Cube whose corners can be correctly oriented, since condition 2 is satisfied, see \cite{Bande82} for details. 

We have proved that the initial configuration $ (id_{S_{8}}, id_{S_{24}}, id_{S_{24}} , 0,0) $ is in the orbit of $ (\sigma,\tau, \rho, x, y) $, hence the latter is valid, concluding the proof of the theorem.

\bibliographystyle{plain}

\end{document}